\documentclass[preprint,3p,12pt]{elsarticle}

\usepackage{lineno,hyperref}
\modulolinenumbers[5]

\usepackage{hyperref}

\usepackage{epstopdf}

\usepackage{mathrsfs}
\usepackage{amssymb}
\usepackage{amsfonts}
\usepackage{latexsym}
\usepackage{amsmath}
\usepackage{xcolor}
\usepackage{wrapfig}
\usepackage{floatflt}

\usepackage{graphicx}
\def\lb{\label}

\newcommand{\er}[1]{\textrm{(\ref{#1})}}

\newtheorem{theorem}{\bf Theorem}[section]
\newtheorem{lemma}[theorem]{\bf Lemma}

\def\a{\alpha}  \def\cA{{\mathcal A}}       
\def\b{\beta}   \def\cB{{\mathcal B}}       
\def\g{\gamma}  \def\cC{{\mathcal C}}       
  \def\cD{{\mathcal D}}       
\def\d{\delta}  \def\cE{{\mathcal E}}

         \def\mL{{\mathscr L}}
\def\l{\lambda}        
        
\def\m{\mu}

\def\s{\sigma}         
         
\def\t{\tau}

       \def\C{{\mathbb C}}    
    \def\N{{\mathbb N}}   

\def\lt{\biggl}                  \def\rt{\biggr}
\def\ol{\overline}               
\def\no{\noindent}

\def\iy{\infty}
\def\sm{\setminus}               \def\es{\emptyset}
\def\ss{\subset}                 \def\ts{\times}
                \def\os{\oplus}

\def\el2{\ell^{\,2}}             \def\1{1\!\!1}

\def\BBox{\hspace{1mm}\vrule height6pt width5.5pt depth0pt \hspace{6pt}}

\newtheorem{corollary}[theorem]{\bf Corollary}

\newcommand{\ca}{\begin{cases}}
\newcommand{\ac}{\end{cases}}
\newcommand{\ma}{\begin{pmatrix}}
\newcommand{\am}{\end{pmatrix}}
\def\eq{\begin{equation}}
\def\qe{\end{equation}}
\def\[{\begin{equation}}
\def\]{\end{equation}}

\def\BBox{\hspace{1mm}\vrule height6pt width5.5pt depth0pt \hspace{6pt}}

\begin{document}

\begin{frontmatter}

\title{Mixed multidimensional integral operators
with piecewise constant kernels and their representations}
\date{\today}

\author
{Anton A. Kutsenko}

\address{Jacobs University (International University Bremen), 28759 Bremen, Germany; email: akucenko@gmail.com}

\begin{abstract}
We consider the algebra of mixed multidimensional integral
operators. In particular, Fredholm integral operators of the first
and second kind belongs to this algebra. For the piecewise constant
kernels we provide an explicit representation of the algebra as a
product of simple matrix algebras. This representation allows us to
compute the inverse operators (or to solve the corresponding
integral equations) and to find the spectrum explicitly. Moreover,
explicit traces and determinants are also constructed. So, roughly
speaking, the analysis of integral operators is reduced to the
analysis of matrices. All the qualitative characteristics of the
spectrum are preserved since only the kernels are approximated.
\end{abstract}

\begin{keyword}
operator algebras, traces and determinants, periodic lattice with
defects, guided and localized waves
\end{keyword}


\end{frontmatter}


{\section{Introduction}\lb{sec1}}

Let me describe briefly the results of the paper. The classical 1D
Fredholm operators
\[\lb{f1}
 \cA u(k)=\int_0^1A(k,x)u(x)dx
\]
are well studied, see, e.g., the famous pioneer work
\cite{Fredholm}. In particular, if
\[\lb{f2}
 A(k,x)=\sum a_{ij}p\chi_i(k)\chi_j(x)
\]
is a piecewise constant ($p$-step, see \er{005}) kernel with
$a_{ij}\in\C$ then the integral operators \er{f1} form an algebra
which is isomorphic to the matrix algebra $\C^{p\ts p}$, see, e.g.,
\cite{Gbook}. The isomorphism is
$\cA\leftrightarrow(a_{ij})_{i,j=1}^p$. But because the algebra of
integral operators does not contain the standard identity operator,
the invertibility in $\C^{p\ts p}$ does not mean the invertibility
in the large algebra of bounded operators. The simplest extension
that devoid of this shortcoming consists of operators of the form
\[\lb{f3}
 \cA u(k)=a u(k)+\int_0^1A(k,x)u(x)dx,\ \ a\in\C.
\]
In this case, the corresponding algebra is isomorphic to
$\C\ts\C^{p\ts p}$, the isomorphism is
$\cA\leftrightarrow(a,(a_{ij}+a\d_{ij})_{i,j=1}^p)$ ($\d$ is the
Kronecker $\d$). The invertibility in $\C\ts\C^{p\ts p}$ means also
the invertibility in the algebra of bounded operators. The result
becomes more substantial if $a$ in \er{f3} is a $p$-step function.
In the current paper, we consider the algebra of all mixed integral
operators of different dimensions with $p$-step kernels that extend
the algebras mentioned above. As it is known, various spectral
problems in multidimensional case are much more complex than in
one-dimensional case. The main motivation for the paper is to show
that discrete analogues of multidimensional integral operators admit
an explicit spectral analysis. We give an explicit representation of
the algebra as a direct product of simple matrix algebras. This
representation leads to explicit procedures of finding inverse
operators and the spectra, and allows us to construct explicit
multidimensional traces and determinants. The algebra of mixed
multidimensional integral operators has many physical applications.
In particular, the $p$-step approximations of periodic operators
acting on the structures with crossing defects of various dimensions
belong to this algebra, see, e.g. \cite{K3,Kjmaa,Kjmaa1} for
operators with defects and \cite{B} for classical periodic operators
without defects.
So, the studying of the guided, surface, and other Rayleigh waves
propagating in such structures are based on the corresponding
determinants and inverse operators. Note that, due to the
Stone-Weierstrass theorem, only step functions allow us to use all
the power of the theory of finite-dimensional algebras. At the same
time, $p$-step kernels can approximate continuous (or other class)
kernels with arbitrary precision when $p\to\iy$. All this shows that
$p$-step kernels are noteworthy, see also \cite{ELKVT}.

Before introducing the algebra of multidimensional integral
operators with step kernels let us define some auxiliary objects. We
fix some $p\in\N$. Let $\a$ be some subset of integer numbers. The
set of multi-indices $P_{\a}$ is defined by
\[\lb{001}
 P_{\a}=\{{\bf m}=(m_n)_{n\in\a}:\ m_n\in\{1,...,p\}\},
\]
where the components of ${\bf m}$ are arranged in the order of
increasing indices $n$. Further, we always assume that the
components of any vectors are arranged in the order of increasing
indices. Let $\a$, $\b$ be two disjoint sets consisting of integer
numbers and let ${\bf x}_{\a}=(x_n)_{n\in\a}$, ${\bf
y}_{\b}=(y_n)_{n\in\b}$ be two vectors. It is convenient to use the
following notation
\[\lb{002}
 {\bf y}_{\b}\diamond{\bf x}_{\a}={\bf x}_{\a}\diamond{\bf y}_{\b}={\bf z}_{\a\cup\b},\ \ where\ \
 z_n=\ca x_n,& n\in\a,\\ y_n,& n\in\b.\ac
\]
Also, for ${\bf x}=(x_n)_{n\in\g}$ and some sets $\a\ss\g$ we denote
${\bf x}_{\a}=(x_n)_{n\in\a}$.
We fix a positive
integer $N$ and denote the set $\eta=\{1,..,N\}$. Let $L^2$ be the
Hilbert space of square-integrable scalar functions acting on the
cube $[0,1)^N$. We consider the operators $\cA: L^2\to L^2$ of the
form
\[\lb{003}
 \cA u({\bf k})=\sum_{\a\ss\eta}\int_{[0,1)^{|\a|}}A_{\a}({\bf k},{\bf
 x}_{\a})u({\bf k}_{\ol{\a}}\diamond{\bf x}_{\a})d{\bf x}_{\a},\ \
 u\in L^2,
\]
where $\ol{\a}=\eta\sm\a$ denotes the complement to the set $\a$,
the vector ${\bf x}=(x_n)_{n\in\eta}$, the differential $d{\bf
x}_{\a}=\prod_{n\in\a}dx_n$ corresponds to the Lebesgue measure, and
$|\a|$ is the number of elements in $\a$. If $\a=\es$ is the empty
set then the corresponding term in \er{003} is $A_{\es}({\bf
k})u({\bf k})$. The $p$-step kernels $A_{\a}$ have the following
form
\[\lb{004}
 A_{\a}({\bf k},{\bf
 x}_{\a})=\sum_{({\bf i},{\bf m})\in P_{\eta}\ts P_{\a}}p^{|\a|}a_{\a}({\bf i},{\bf
 m})\prod_{j\in\eta}\chi_{i_j}(k_j)\prod_{n\in\a}\chi_{m_n}(x_n),
\]
where $a_{\a}({\bf i},{\bf m})\in\C$ and the step functions $\chi$
are defined by
\[\lb{005}
 \chi_r(y)=\ca 1,& y\in[(r-1)/p,r/p),\\ 0, & otherwise. \ac
\]
The operators $\cA$ \er{003} form an algebra $\mL$ which is a
subalgebra of the algebra of bounded operators acting in $L^2$. In
fact, this algebra is generated by the following $pN+p$ elementary
operators
\[\lb{006}
 \mL={\rm Alg}\lt(\lt\{\chi_i(k_j)\cdot,\int_0^1\cdot
 dx_j\rt\}_{i=1,j=1}^{p,N}\rt),
\]
where $\cdot$ denotes the place of the operator argument $u\in L^2$.
Note that $\mL$ contains the identical operator from the large
algebra of bounded operators. Hence, the invertibility in $\mL$
means also the invertibility in the algebra of bounded operators and
vice versa. Because $\mL$ is a finite dimensional Von Neumann
algebra, it can be represented as a direct product of simple
algebras. Our goal is to find this representation explicitly. There
are probably various ways to do that, one based on composition
series is preferred for us since it allows us to simplify some of
calculations. We denote the simple matrix algebras as $\C^{n\ts n}$,
$n\in\N$ ($\C^{1\ts1}=\C$). For any $\a\ss\eta$, ${\bf i}\in
P_{\eta}$, ${\bf m}\in P_{\a}$ introduce the following functions
\[\lb{007}
 b_{\a}({\bf i},{\bf
 m})=\sum_{\b\ss\a}\d({\bf i}_{\a\sm\b},{\bf m}_{\a\sm\b})a_{\b}({\bf i},{\bf m}_{\b}),
\]
where the Kronecker delta satisfies: $\d(x,y)=1$ if $x=y$ and
$\d(x,y)=0$ if $x\ne y$. Next, for any $\a\ss\eta$ and ${\bf i}\in
P_{\ol{\a}}$ introduce the matrices
\[\lb{008}
 {\bf B}_{\a}({\bf i})=(b_{\a}({\bf i}\diamond{\bf m},{\bf n}))_{{\bf m},{\bf n}\in
 P_{\a}}\in\C^{p^{|\a|}\ts p^{|\a|}}.
\]
Introduce the following mapping
\[\lb{009}
 \s:\mL\to\prod_{n=0}^N(\C^{p^n\ts p^n})^{{\binom {N} {n}}p^{N-n}},\ \
 \s(\cA)=(({\bf B}_{\a}({\bf i}))_{{\bf i}\in P_{\ol{\a}}})_{\a\ss\eta},
\]
where $\cA$, ${\bf B}$ are defined in \er{003}-\er{004},
\er{007}-\er{008} and $\binom {N} {n}$ are the binomial
coefficients. The identities \er{007}-\er{009} allow us to compute
the inverse mapping $\s^{-1}$ explicitly by using
\[\lb{inverse}
 a_{\a}({\bf i},{\bf m})=\sum_{\b\ss\a}(-1)^{|\a\sm\b|}\d({\bf
 i}_{\a\sm\b},{\bf m}_{\a\sm\b})b_{\b}({\bf i},{\bf m}_{\b}).
\]
The following theorem is our main result.

\begin{theorem}\lb{T1}
The mapping $\s$  is an algebra isomorphism.
\end{theorem}

We immediately obtain the next corollary which is useful in physical
and numerical applications.

\begin{corollary}\lb{C1}
i) The operator $\cA$ is invertible if and only if all matrices
${\bf B}_{\a}({\bf i})$ are invertible. In this case
\[\lb{010}
 \cA^{-1}=\s^{-1}(({\bf B}_{\a}^{-1}({\bf i}))_{{\bf i}\in P_{\ol{\a}}})_{\a\ss\eta}.
\]
ii) The spectrum of $\cA$ consists of all eigenvalues of the
matrices ${\bf B}_{\a}({\bf i})$. The eigenvalues of ${\bf
B}_{\eta}$ form a discrete spectrum, other eigenvalues belong to
essential spectrum.

\no iii) The multidimensional trace $\t$ and the determinant $\pi$
can be defined as follows
\[\lb{011}
 \t,\pi:\mL\to\prod_{n=0}^N\C^{{\binom {N} {n}}p^{N-n}}=\C^{(p+1)^N},
\]
\[\lb{012}
 \t(\cA)=(({\rm Tr}{\bf B}_{\a}({\bf i}))_{{\bf i}\in P_{\ol{\a}}})_{\a\ss\eta},\ \
 \pi(\cA)=(({\rm det}{\bf B}_{\a}({\bf i}))_{{\bf i}\in P_{\ol{\a}}})_{\a\ss\eta}.
\]
They satisfy the usual properties
\[\lb{013}
 \t(\l\cA+\m\cB)=\l\t(\cA)+\m\t(\cB),\ \ \t(\cA\cB)=\t(\cB\cA),\ \
 \pi(\cA\cB)=\pi(\cA)\pi(\cB),\ \ \pi(e^{\cA})=e^{\t(\cA)},
\]
where $\l,\m\in\C$ and $\cA,\cB\in\mL$. The operator is invertible
iff all determinants are nonzero.
\end{corollary}

{\bf Remark.} All the results can be easily extended to the case
where $\cA$ \er{003} acts on $\os_{m=1}^ML^2$ and the kernels ${\bf
A}_{\a}$ \er{004} are $M\ts M$ $p$-step matrices (this means that
all entries of these matrices are $p$-step functions and, hence, the
coefficients ${\bf a}_{\a}$ \er{004} are $M\ts M$ constant
matrices). The integral operators with $p$-step $M\ts M$
matrix-valued kernels ${\bf A}_{\a}$ form an algebra $\mL_M$. It can
be shown that
\[\lb{mult}
 \mL_M \simeq \prod_{n=0}^N(\C^{Mp^n\ts Mp^n})^{{\binom {N} {n}}p^{N-n}}
\]
and the corresponding isomorphism has the same form as $\s$ \er{009}
but with ${\bf a}_{\a}$ instead of $a_{\a}$.

\section{\lb{S1}Proof of Theorem \ref{T1}}
Introduce the following operators
\[\lb{100}
 \cA_{ij}=\chi_{i}(k_j)\cdot,\ \ \cB_{ij}=p\int_0^1\chi_i(x_j)\cdot
 dx_n.
\]

\begin{lemma}\lb{L1} The operators $\cA,\cB$ are commute if their
second indices are different. Moreover, for any $i,j,r$ the
following identities hold true
\[\lb{101}
 \cA_{ij}\cA_{rj}=\d(i,r)\cA_{ij},\ \
 \cB_{ij}\cA_{rj}=\d(i,r)\cB_{ij},\ \ \cB_{ij}\cB_{rj}=\cB_{rj}.
\]
\end{lemma}
{\it Proof.} The direct calculations and Fubini theorem give the
result. \BBox

For any $\a\ss\eta$, ${\bf i}\in P_{\eta}$, ${\bf m}\in P_{\a}$
introduce the following operators
\[\lb{102}
 \cC_{\a}({\bf i},{\bf
 m})=\prod_{j\in\eta}\cA_{i_jj}\prod_{j\in\a}\cB_{m_jj},\
 \ \cD_{\a}({\bf i},{\bf m})=\cC_{\a}({\bf i},{\bf
 m})\prod_{j\in\ol{\a}}(1-\cB_{i_jj}).
\]
These operators are connecting through the next equations.
\begin{lemma}\lb{L2}
For any $\a$, ${\bf i}$, ${\bf m}$ the following identities hold
true
\[\lb{ident}
 \cC_{\a}({\bf i},{\bf m})=\sum_{\b\supset\a}\cD_{\b}({\bf i},{\bf m}\diamond{\bf
 i}_{\b\sm\a}),\ \ \ \cD_{\a}({\bf i},{\bf m})=
 \sum_{\b\supset\a}(-1)^{|\a\sm\b|}\cC_{\b}({\bf i},{\bf m}\diamond{\bf
 i}_{\b\sm\a}).
\]
\end{lemma}
{\it Proof.} Using the commutativity from Lemma \ref{L1} and
\er{102} we deduce that
\[\lb{i1}
 \cC_{\a}({\bf i},{\bf m})=\cC_{\a}({\bf i},{\bf
 m})\prod_{j\in\ol{\a}}(1-\cB_{i_jj}+\cB_{i_jj})=
 \sum_{\b\supset\a}\cD_{\b}({\bf i},{\bf m}\diamond{\bf
 i}_{\b\sm\a})
\]
and
\[\lb{i2}
 \cD_{\a}({\bf i},{\bf m})=\cC_{\a}({\bf i},{\bf
 m})\prod_{j\in\ol{\a}}(1-\cB_{i_jj})=\sum_{\b\supset\a}(-1)^{|\a\sm\b|}\cC_{\b}({\bf i},{\bf m}\diamond{\bf
 i}_{\b\sm\a}).\ \BBox
\]

While $\cC_{\a}({\bf i},{\bf m})$ is a standard basis in $\mL$ (see
\er{003}-\er{004}), the basis $\cD_{\a}({\bf i},{\bf m})$ is
"orthogonal" that will be proved in the next lemma.

\begin{lemma}\lb{L3}
The following identities hold true
\[\lb{103}
 \cD_{\a}({\bf i},{\bf m})\cD_{\b}({\bf p},{\bf r})=\ca 0, &
 \a\ne\b,\\ 0,& \a=\b,\ {\bf i}_{\ol{\a}}\ne{\bf p}_{\ol{\a}},\\
 \d({\bf m},{\bf p}_{\a})\cD({\bf i},{\bf r}),& \a=\b,\ {\bf i}_{\ol{\a}}={\bf
 p}_{\ol{\a}}. \ac
\]
\end{lemma}
{\it Proof.} The case $\a\ne\b$. There are two possibilities. The
first one, there is $i_j\in\ol{\a}\sm\ol{\b}$. Due to the
commutativity following from Lemma \ref{L1}, we have
\[\lb{104}
 \cD_{\a}({\bf i},{\bf m})\cD_{\b}({\bf p},{\bf
 r})=\cE\cA_{i_jj}(1-\cB_{i_jj})\cA_{p_jj}\cB_{r_jj},
\]
where $\cE$ is a product of elements $\cA$, $\cB$, $1-\cB$ which
have second indices not equal to $j$. Using \er{101} we obtain the
following identities
\[\lb{105}
 \cA_{i_jj}(1-\cB_{i_jj})\cA_{p_jj}\cB_{r_jj}=\d(i_j,p_j)\cA_{i_jj}(1-\cB_{i_jj})\cB_{r_jj}=
 \d(i_j,p_j)\cA_{i_jj}(\cB_{r_jj}-\cB_{r_jj})=0.
\]
The second possibility, there is $p_j\in\ol{\b}\sm\ol{\a}$. Due to
the commutativity following from Lemma \ref{L1}, we have
\[\lb{106}
 \cD_{\a}({\bf i},{\bf m})\cD_{\b}({\bf p},{\bf
 r})=\cE\cA_{i_jj}\cB_{m_jj}\cA_{p_jj}(1-\cB_{p_jj}),
\]
where $\cE$ is a product of elements $\cA$, $\cB$, $1-\cB$ which
have second indices not equal to $j$. Using \er{101} we obtain the
following identities
\[\lb{107}
 \cA_{i_jj}\cB_{m_jj}\cA_{p_jj}(1-\cB_{p_jj})=\d(m_j,p_j)\cA_{i_jj}\cB_{m_jj}(1-\cB_{p_jj})
 =\d(m_j,p_j)\cA_{i_jj}(\cB_{m_jj}-\cB_{p_jj})=0
\]
since $\d(m_j,p_j)=0$ for $m_j\ne p_j$ and $\cB_{m_jj}-\cB_{p_jj}=0$
for $m_j=p_j$. We have proved the first identity in \er{103}.

The case $\a=\b$, ${\bf i}_{\ol{\a}}\ne{\bf p}_{\ol{\a}}$. Then
there is $j\in\ol{\a}=\ol{\b}$ such that $i_j\ne p_j$. Due to the
commutativity (see Lemma \ref{L1}), we have
\[\lb{109}
 \cD_{\a}({\bf i},{\bf m})\cD_{\b}({\bf p},{\bf
 r})=\cE\cA_{i_jj}(1-\cB_{i_jj})\cA_{p_jj}(1-\cB_{p_jj}),
\]
where $\cE$ is a product of elements $\cA$, $\cB$, $1-\cB$ which
have second indices not equal to $j$. Using \er{101} we obtain
\[\lb{110}
 \cA_{i_jj}(1-\cB_{i_jj})\cA_{p_jj}(1-\cB_{p_jj})=\d(i_j,p_j)\cA_{i_jj}(1-\cB_{i_jj})(1-\cB_{p_jj})=0.
\]
We have proved the first identity in \er{103}.

The case $\a=\b$, ${\bf i}_{\ol{\a}}={\bf p}_{\ol{\a}}$. Due to the
commutativity (see Lemma \ref{L1}), we have
\[\lb{111}
 \cD_{\a}({\bf i},{\bf m})\cD_{\b}({\bf p},{\bf
 r})=\prod_{j\in\ol{\a}}(\cA_{i_jj}(1-\cB_{i_jj})\cA_{p_jj}(1-\cB_{p_jj}))
 \prod_{j\in\a}(\cA_{i_jj}\cB_{m_jj}\cA_{p_jj}\cB_{r_jj}).
\]
Using \er{101} and $i_j=p_j$ for $j\in\ol{\a}$ we obtain
\[\lb{112}
 \cA_{i_jj}(1-\cB_{i_jj})\cA_{p_jj}(1-\cB_{p_jj})=\cA_{i_jj}(1-\cB_{i_jj})(1-\cB_{i_jj})
\]
\[\lb{113}
 =\cA_{i_jj}(1-\cB_{i_jj}-\cB_{i_jj}+\cB_{i_jj})=\cA_{i_jj}(1-\cB_{i_jj}).
\]
Using \er{101} we obtain also
\[\lb{114}
 \cA_{i_jj}\cB_{m_jj}\cA_{p_jj}\cB_{r_jj}=\d(m_j,p_j)\cA_{i_jj}\cB_{r_jj}.
\]
Substituting \er{112}-\er{114} into \er{111} leads to
\[\lb{115}
 \cD_{\a}({\bf i},{\bf m})\cD_{\b}({\bf p},{\bf
 r})=\prod_{j\in\a}\d(m_j,p_j)\prod_{j\in\ol{\a}}(\cA_{i_jj}(1-\cB_{i_jj}))\prod_{j\in\a}(\cA_{i_jj}\cB_{r_jj})
\]
\[\lb{116}
 =\d({\bf m},{\bf p}_{\a})\cC({\bf i},{\bf
 r})\prod_{j\in\ol{\a}}(1-\cB_{i_jj})=\d({\bf m},{\bf p}_{\a})\cD({\bf i},{\bf
 r}),
\]
where we have used the commutativity from Lemma \ref{L1}. Now
\er{103} is completely proved. \BBox

{\bf Proof of Theorem \ref{T1}.} The identities \er{103} mean that
\[\lb{117}
 {\rm Alg}(\bigcup_{\a\ss\eta}\bigcup_{{\bf i}\in P_{\eta}}\bigcup_{{\bf m}\in P_{\a}}\{\cD_{\a}({\bf i},{\bf
 m})\})=\prod_{\a\ss\eta}{\rm Alg}(\bigcup_{{\bf i}\in P_{\eta}}\bigcup_{{\bf m}\in P_{\a}}\{\cD_{\a}({\bf i},{\bf
 m})\})
\]
\[\lb{118}
 =\prod_{\a\ss\eta}\prod_{{\bf p}\in P_{\ol{\a}}}{\rm Alg}(\bigcup_{{\bf r}\in P_{\a}}
 \bigcup_{{\bf m}\in P_{\a}}\{\cD_{\a}({\bf p}\diamond{\bf r},{\bf
 m})\})=\prod_{\a\ss\eta}\prod_{{\bf p}\in P_{\ol{\a}}}\C^{p^{|\a|}\ts
 p^{|\a|}}=\prod_{n=0}^N(\C^{p^n\ts p^n})^{{\binom {N} {n}}p^{N-n}},
\]
where $\cD_{\a}({\bf p}\diamond{\bf r},{\bf m})$ corresponds to the
elementary matrix $(\d({\bf r},{\bf m}))_{{\bf r},{\bf m}\in
P_{\a}}$ in the algebra $\C^{|P_{\a}|\ts|P_{\a}|}=\C^{p^{|\a|}\ts
p^{|\a|}}$. Let $\cA\in\mL$ be some operator of the form
\er{003}-\er{004}. Using \er{102}, \er{ident} we deduce that
\[\lb{119}
 \cA=\sum_{\a\ss\eta}\sum_{{\bf i}\in P_{\eta}}\sum_{{\bf m}\in P_{\a}}a_{\a}({\bf i},{\bf m})\cC_{\a}({\bf i},{\bf
 m})=\sum_{\a\ss\eta}\sum_{{\bf i}\in P_{\eta}}\sum_{{\bf m}\in P_{\a}}a_{\a}({\bf i},{\bf
 m})\sum_{\b\supset\a}\cD_{\b}({\bf i},{\bf m}\diamond{\bf
 i}_{\b\sm\a})
\]
\[\lb{120}
 =\sum_{\a\ss\eta}\sum_{{\bf i}\in P_{\eta}}\sum_{{\bf m}\in
 P_{\a}}(\sum_{\b\ss\a}\d({\bf
 i}_{\a\sm\b},{\bf m}_{\a\sm\b})
 a_{\b}({\bf i},{\bf m}_{\b}))\cD_{\a}({\bf i},{\bf m})=
 \sum_{\a\ss\eta}\sum_{{\bf i}\in P_{\eta}}\sum_{{\bf m}\in
 P_{\a}}b_{\a}({\bf i},{\bf m})\cD_{\a}({\bf i},{\bf m})
\]
which give us the form of matrices ${\bf B}_{\a}$ \er{008} in the
basis $\cD_{\a}$. To calculate the inverse mapping $\s^{-1}$ we take
the operator $\cA$ in the basis $\cD_{\a}$ and, using \er{ident},
write it in the standard basis $\cC_{\a}$, i.e.
\[\lb{121}
 \cA=\sum_{\a\ss\eta}\sum_{{\bf i}\in P_{\eta}}\sum_{{\bf m}\in
 P_{\a}}b_{\a}({\bf i},{\bf m})\cD_{\a}({\bf i},{\bf m})
 =\sum_{\a\ss\eta}\sum_{{\bf i}\in P_{\eta}}\sum_{{\bf m}\in
 P_{\a}}b_{\a}({\bf i},{\bf m})(\sum_{\b\supset\a}(-1)^{|\a\sm\b|}\cC_{\b}({\bf i},{\bf m}\diamond{\bf
 i}_{\b\sm\a}))
\]
\[\lb{122}
 =\sum_{\a\ss\eta}\sum_{{\bf i}\in P_{\eta}}\sum_{{\bf m}\in
 P_{\a}}(\sum_{\b\ss\a}(-1)^{|\a\sm\b|}\d({\bf
 i}_{\a\sm\b},{\bf m}_{\a\sm\b})b_{\b}({\bf i},{\bf m}_{\b}))\cC_{\a}({\bf i},{\bf m})
 =\sum_{\a\ss\eta}\sum_{{\bf i}\in P_{\eta}}\sum_{{\bf m}\in
 P_{\a}}a_{\a}({\bf i},{\bf m})\cC_{\a}({\bf i},{\bf m})
\]
which give us \er{inverse}. \BBox

\section{\lb{S2}Examples}

We consider the simplest case $p=1$
\[\lb{200}
 \mL={\rm Alg}\lt(1,\int_0^1\cdot dx_1,...,\int_0^1\cdot
 dx_N\rt)\simeq\C^{2^N}.
\]
By Corollary \ref{C1} the operator
\[\lb{201}
 \cA u({\bf k})=\sum_{\a\ss\{1,...,N\}}a_{\a}\int_{[0,1)^{|\a|}}u({\bf k}) d{\bf
 k}_{\a},\ \ u\in L^2_N
\]
is invertible if and only if the numbers
$
 b_{\a}=\sum_{\b\ss\a}a_{\b}
$
are all non-zero and then (see \er{inverse})
\[\lb{203}
 \cA^{-1} u({\bf k})=
 \sum_{\a\ss\{1,...,N\}}\lt(\sum_{\b\ss\a}(-1)^{|\a\sm\b|}b^{-1}_{\b}\rt)\int_{[0,1)^{|\a|}}u({\bf k}) d{\bf
 k}_{\a},\ \ u\in L^2_N.
\]
In particular, for $N=2$ we have
\[\lb{204}
 \lt(a\cdot+b\int_0^1\cdot dk_1+c\int_0^1\cdot dk_2+d\int_0^1\int_0^1\cdot
 dk_1dk_2\rt)^{-1}=a^{-1}\cdot-
\]
\[\lb{205}
 \frac{b}{a(a+b)}\int_0^1\cdot dk_1-\frac{c}{a(a+c)}\int_0^1\cdot
 dk_2+\frac{(2a+b+c+d)bc-a^2d}{a(a+b)(a+c)(a+b+c+d)}\int_0^1\int_0^1\cdot
 dk_1dk_2.
\]
Identity \er{203} leads also to the following example
\[\lb{206}
 \lt(\sum_{\a\ss\{1,...,N\}}\int_{[0,1)^{|\a|}}\cdot d{\bf
 k}_{\a}\rt)^{-1}=\sum_{\a\ss\{1,...,N\}}(-1/2)^{|\a|}\int_{[0,1)^{|\a|}}\cdot d{\bf
 k}_{\a}.
\]

\section*{Acknowledgements}
This work was partially supported by the RSF project
N\textsuperscript{\underline{o}}15-11-30007 and TRR 181 project.

\bibliography{bibl_perp1}

\begin{thebibliography}{1}
\expandafter\ifx\csname url\endcsname\relax
  \def\url#1{\texttt{#1}}\fi
\expandafter\ifx\csname urlprefix\endcsname\relax\def\urlprefix{URL }\fi
\expandafter\ifx\csname href\endcsname\relax
  \def\href#1#2{#2} \def\path#1{#1}\fi

\bibitem{Fredholm}
{E.~I.~Fredholm}, {Sur une classe d'equations fonctionnelles}, Acta Math. 27
  (1903) 365--390.

\bibitem{Gbook}
M.~S. Gockenbach, Finite-Dimensional Linear Algebra, Discrete mathematics and
  its applications., CRC Press, 2010.

\bibitem{K3}
A.~A. Kutsenko, Analytic formula for amplitudes of waves in lattices with
  defects and sources and its application for defects detection, Eur. J. Mech.
  A-Solid. 54 (2015) 209--217.

\bibitem{Kjmaa}
{A.~A.~Kutsenko}, Algebra of multidimensional periodic operators with defects,
  J. Math. Anal. Appl. 428 (2015) 221--230.

\bibitem{Kjmaa1}
{A.~A.~Kutsenko}, Algebra of 2d periodic operators with local and perpendicular
  defects, J. Math. Anal. Appl. 442 (2016) 796--803.

\bibitem{B}
{L.~Brillouin}, Wave propagation in periodic structures, {Dover Publications
  Inc, New York}, 2003.

\bibitem{ELKVT}
{E.~Espinoza-Loyola, Yu.~I.~Karlovich, O.~Vilchis-Torres}, {C$^*$-algebras of
  Bergman type operators with piecewise constant coefficients over sectors},
  Integr. Equ. Oper. Theory 83 (2015) 243--269.

\end{thebibliography}

\end{document}